%

\hsize=6in
\vsize=7.5in
\hoffset=.125in
\font\Bbb=msbm10 at 12pt
\font\ger=eufm10 at 12pt
\font\Large=cmbx10 at 15pt 
\parskip=.2in
\baselineskip=10pt
\magnification=1200
{\fiverm This is the Plain-Tex version of a paper that has been accepted
for publication in the} {\fivebf Journal of Natural Geometry}{\fiverm --an
international quarterly that started publication in 1992 with the stated aim
of bringing unity through simplicity in the mathematics of natural structures.
The address of the journal can be found in the indices of the Mathematical
Reviews.}\par

\baselineskip=18pt

\parskip=.2in




\baselineskip=24pt
{\centerline {\Large  On a Combinatorial Property of Families of}} 
{\centerline {\Large of  Sequences  Converging to $\infty$}}

\baselineskip=18pt
\vskip 32pt
{\centerline {\Large  Apoloniusz Tyszka }}

\baselineskip=18pt

{\parindent=.5in \narrower
\noindent  {\bf Abstract.} We consider  families $\Phi $  of  sequences
converging to $\infty $, with the property that  for  every  open   
set $U\subseteq${\Bbb R} that is unbounded above there exists a sequence 
belonging to $\Phi $,
which has an infinite number of terms belonging to $U$.

\par}

\baselineskip=16pt


In this paper we consider  families $\Phi $ of sequences converging to 
$\infty $ that $\Phi $ satisfy the following condition (C):
\par
\noindent (C):\quad  {\it for every open set} 
$U\subseteq$ {\Bbb R} {\it that is unbounded above there exists  a  sequence
belonging to} $\Phi ${\it , which has an infinite number of terms 
belonging to U.}
\par
\noindent The origin of this consideration  is  the observation that  
if  a  family $\Phi $  of
sequences converging to $\infty $  satisfies  
the  condition  (C)  and $f:$ {\Bbb R}$\rightarrow${\Bbb R}  is
continuous,  then  convergence  of $\{f(x_{n})\}_{n\in \omega }$  
to $d~\in ~${\Bbb R}$\cup \{-\infty ,\infty \}$  for   all
sequences $\{x_{n}\}_{n\in \omega }$ belonging to $\Phi $ implies 
that $\lim_{x\rightarrow \infty }f(x)$~=~{\it d}.
\par
We present below necessary definitions from [1] (see also  [2],  [3]).
For the functions $f,g:\omega \rightarrow \omega $ we define:
$$ f~\le_*g \Leftrightarrow ~\{i:~f(i)~>~g(i)\}\;{\rm is\; finite\;}
(g\;{\rm dominates\;}f).$$
\noindent We say that $F~\subseteq \omega^\omega $ is  bounded  
if $\exists ~g:\omega \rightarrow \omega ~\forall ~f\in F~f~\le_*g$ 
(briefly: $F~\le_*g$,
$g$~dominates $F$). Every countable set $F=\{f_j:j\in \omega \}$  
is  bounded,  because  the
function $g(n)~=~\max ~\{f_{j}(n):~j\le n\}~+~1$ dominates 
$F$ (i.e. $F~\le_* g$).
Let ${\bf b~}=~\min ~\{$card~$F: F$ is unbounded$\}$.
\par
%
\noindent{\bf Theorem 1.} {\it If} $\Phi $  {\it is  a  family  of  sequences  
converging  to} $\infty $  {\it and}
card~$\Phi ~<~{\bf b}${\it , then} $\Phi $ {\it does not satisfy  
condition} (C){\it .}\hfill \break
{\it Proof.} Let card~$\Phi ~=~\kappa ,
\; \Phi ~=~\{\phi_\alpha:~\alpha <\kappa \}$. There  exists  such  a  
sequence $\{a_n\}_{n\in \omega }$ strictly increasing and converging to 
$\infty $, that for every $\alpha <\kappa $ the set
$\{a_n:~n\in \omega \}$ is disjoint from the set of terms of the sequence 
$\phi _{\alpha }$.  Let $h_{\alpha }(n)$
be the least positive integer for which 
$\left(a_n-~{1/{h_{\alpha}(n)}},\;~a_n+~{1/~ 
{h_{\alpha }(n)}}\right)$ is disjoint from the set of terms of the sequence 
$\phi _{\alpha }$. Since  the  family $\{h_\alpha :~\alpha <\kappa \}$  is
bounded, there exists such a function $h:\omega \rightarrow \omega $ that 
for every $\alpha <\kappa\;  h_{\alpha }\le_{*} h$.  In
the open  set $\bigcup_{n\in \omega} 
\left(a_n-~{1/h(n)}~,\;~a_n+~{1/~ 
h(n)}\right)$
 every sequence $\phi _{\alpha }$
has only a finite number of terms.\hfill $\bullet$
\par
\noindent{\bf Corollary.} {\it Every family of sequences converging to} 
$\infty $  {\it which  satisfies
the condition} (C) {\it is uncountable, Martin's axiom  implies  that  each  
such family  has  cardinality} c ({\it because  Martin's  axiom  implies  
that} {\bf b}$=${\ger c}{\it , see~}[1]{\it ,~}[3]).
\par
For a sequence $\{a_n\}_{n\in \omega }$converging to $\infty $ we define the 
non-decreasing function \hfill\break
$f_{\{a_n\}}:~\omega \rightarrow \omega $ : if $\bigcup_{n\in \omega }
(a_{n}-1,a_{n}+1)~\supseteq ~(i,\infty )$  then
$$
f_{\{a_n\}}(i)~=~\max \left\{ j\in \omega \backslash \{0\}~:\bigcup_{n \in \omega}
\left(a_n-~{1\over j},\;~a_n+~{1\over j}\right)
\supseteq ~(i,\infty )\right\}$$
else $f_{\{a_n\}}(i)~=~0$.
The following Lemma is obvious.
\par
\noindent{\bf Lemma.} {\it If for a family} $\Phi $ {\it of sequences 
converging to} $\infty $  {\it the  family  of
functions}\hfill\break 
$\left\{ f_{\{a_{n}\}}~:~\{a_{n}\}\in \Phi\right\} $   
{\it is   unbounded,   then} $\Phi $   {\it satisfies condition~}(C){\it .}
\par
\noindent {\bf Theorem 2.} {\it There exists a family of sequences 
converging  to} $\infty $  {\it which
satisfies  condition} (C) {\it and has cardinality} {\bf b}.\hfill \break
{\it Proof.} Let $F~\subseteq {\omega}^\omega $ be an unbounded  
family  and  card~$F$~=~{\bf b}.  Replacing
every function $f\in F$ by the function
$\tilde{f}(n)=\max \{f(0),~f(1),...,f(n)\}$ we  obtain
an unbounded family $\tilde{F}$  of non-decreasing functions.  
Hence,  card$\tilde{F}~=~${\bf b}.  To
every function $g\in \tilde{F}$  we assign a sequence converging to $\infty $:
\par
\noindent $0,~~{1\over g(0)+1},~~{2\over g(0)+1},~~...,~~{g(0)\over g(0)+1},~~
1,~~1{1\over g(1)+1},~~1{2\over g(1)+1},~~ ...,~~1{g(1)\over g(1)+1},~~
2,~~2{1\over g(2)+1}$,\par

\noindent$2{2\over g(2)+1},~~...,~~2{g(2)\over g(2)+1}${\it ,...}
\par
\noindent We obtain the family of sequences which has cardinality {\bf b}, according to the
Lemma this family satisfies the condition (C).\hfill$\bullet$
\par
{\bf Theorem 3.}  {\it If} $X~\subseteq ~${\Bbb R} {\it is a set of second 
category, then the  family  of
sequences} $\{\{x+\log (n+1)\}_{n \in \omega}:x\in X\}$ 
{\it satisfies  condition} (C).\hfill \break
\noindent {\it Proof.} Transformation $x\rightarrow \exp (-x)$ allows us 
to formulate  this  theorem
in an equivalent form:
\par
\noindent {\it If $X~\subseteq ~(0,\infty )$ is a set of second category, then for every open set $U~\subseteq ~(0,\infty )$
with $0\in \overline{U}$  there exists an $x\in X$ such that an 
infinite number of terms  of  the
sequence $\{{x\over n}\}_{n\in \omega \backslash\{ 0 \}}$ belong to }$U$.
\par
 Using operations on sets we can formulate this condition as follows:
$$\left[\bigcap^{}_{n\in \omega \backslash \{0\}} \bigcup^{}_{k\ge n}k
U \right]\cap ~X~\neq \emptyset.$$

Every set $\bigcup^{}_{k\ge n}kU$ is an open and dense subset of 
$[0,\infty )$,  so  we  establish  the
asertion using the Baire category theorem.\hfill$\bullet$
\par
The origin of Theorem 3 is the problem number  107  from  the  problem
book [5]:
\par
{\sevenrm
\noindent {\it Assume that} $f:(0,\infty )\rightarrow$ {\Bbb R} 
{\it is continuous and for every} $x>0$ {\it the limit} 
$\lim_{n\rightarrow \infty }f\left({x\over n}\right)$
{\it exists and equals 0. Is it possible to conclude that} 
$\lim_{n\rightarrow 0^+}f(x)~=~0${\it ?}}
\par
\noindent From Theorem 3 follows a positive solution of this problem under  
the  weaker
assumption that the set 
$\left\{ x \in (0,\infty ):~\lim_{n\rightarrow \infty }f\left({x\over n}~
\right)=0\right\}$  is a set of second category.
\par
\noindent{\bf Remark} [4]. If a sequence $\{a_n\}_{n\in \omega }$ converging to $\infty $ satisfies
$$\lim_{n\rightarrow \infty}(a_{n+1}-a_{n})~=~0 \eqno (1)$$
 then
\par
for every open  set $U~\subseteq ~${\Bbb R} that is unbounded above the set\hfill \break
$\phantom{xx.}\{r\in ${\Bbb R}$\;:$ card$~\{r+a_{n}:n\in \omega \}\cap U~\ge ~\omega \}$ 
is a dense subset of {\Bbb R};\hfill (2)\par
\noindent if we additionally assume that  a  sequence 
$\{a_n\}_{n\in \omega }$  converging  to $\infty $  is
non-\hfill \break decreasing then (2) implies (1).
\par
The author  suggests  that  if  a  sequence $\{a_n\}_{n\in \omega }$  
converging  to $\infty $
satisfies (1) then a  dense  subset  considering  in  (2)  is  a  countable
intersection of open dense sets, so if $X~\subseteq ~${\Bbb R} 
is a set of  second  category,
then the family of sequences $\{\{x+a_{n}\}_{n \in \omega}:x\in X\}$ 
satisfies the condition (C).
\par
Using the Baire category theorem we can prove that if $X\subseteq$ {\Bbb R}
is a set of the second category and if a continuous function $h:${\Bbb R}$\rightarrow
[0,1]$ is periodic and surjective, then the family of sequences 
$\left\{\{n+h(nx)\}_{n\in \omega}:x\in X\right\}$ satisfies condition (C).
Using the Baire category theorem for the space of irrationals we can replace
$h(nx)$ by the fractional part of $nx$ in the statement above.
\medskip

\noindent{\Large References}


\item{{\bf 1.}}  R.~Frankiewicz and P.~Zbierski,  {\it Hausdorff gaps and limits,}
\hfill\break Amsterdam: North-Holland, 1994.

\item{{\bf 2.}}  S.~Hechler, On the existence  of  certain  cofinal  subsets  of $^{\omega }\omega $,
In {\it Axiomatic Set Theory,} (ed: T. Jech),  Proc. Symp. Pure Math.{ \bf 13}, II (1974) 155-173
(Amer. Math. Soc.).

\item{{\bf 3.}} T.~Jech, {\it Set theory,} New York: Academic Press, 1978.

\item{{\bf 4.}} Z. Marciniak, {\it Private communication,} 1991.

\item{{\bf 5.}}  D.~J.~Newman,  {\it A  Problem  Seminar}, Heidelberg:  Springer,1982.

\noindent {\bf DEPARTMENT OF AGRICULTURAL MECHANIZATION,\hfill \break AGRICULTURAL UNIVERSITY,
BALICKA 104, 30-149 KRAK\'OW,\hfill \break POLAND.}
{\it email: rttyszka@pl.edu.cyf-kr}

\hfill {\sevenbf RECEIVED BY THE EDITOR:  NOVEMBER 15, 1994}
\end